\documentclass[letterpaper, 11pt]{amsart}
\usepackage[all]{xy}        
\CompileMatrices            
\UseTips                    

\usepackage[german,english]{babel}   
\usepackage[bookmarks=true]{hyperref}   
\usepackage{amssymb}
\usepackage{xspace}

\theoremstyle{plain}
\newtheorem{theorem}{Theorem}[section]
\newtheorem*{theorem*}{Theorem}
\newtheorem{proposition}[theorem]{Proposition}

\newtheorem{lemma}[theorem]{Lemma}

\theoremstyle{definition}
\newtheorem{definition}[theorem]{Definition}

\theoremstyle{remark}
\newtheorem{remark}[theorem]{Remark}
\newtheorem{example}[theorem]{Example}

\numberwithin{equation}{section}

\newcommand{\enm}[1]{\ensuremath{#1}}          
\newcommand{\op}[1]{\operatorname{#1}}
\newcommand{\cal}[1]{\mathcal{#1}}

\newcommand{\CC}{\enm{\mathbb{C}}}             

\newcommand{\ZZ}{\enm{\mathbb{Z}}}

\newcommand{\LL}{\enm{\mathbb{L}}}

\newcommand{\Cc}{\enm{\cal{C}}}
\newcommand{\Dd}{\enm{\cal{D}}}

\newcommand{\Hh}{\enm{\cal{H}}}

\newcommand{\Ll}{\enm{\cal{L}}}
\newcommand{\Mm}{\enm{\cal{M}}}

\newcommand{\Oo}{\enm{\cal{O}}}

\renewcommand{\phi}{\varphi}        
\renewcommand{\theta}{\vartheta}
\renewcommand{\epsilon}{\varepsilon}

\newcommand{\Spec}{\op{Spec}}

\newcommand{\Hom}{\op{Hom}}

\newcommand{\Image}{\op{Im}}

\newcommand{\tensor}{\otimes}         
\newcommand{\dirsum}{\oplus}

\renewcommand{\to}[1][]{\xrightarrow{\ #1\ }}

\newcommand{\into}{\hookrightarrow}

\newcommand{\usc}[1][m]{\underline{\phantom{#1}}}

\newcommand{\defeq}{\stackrel{\scriptscriptstyle \op{def}}{=}}

\newcommand{\ie}{\textit{i.e.}\ }           

\newcommand{\Sol}{\op{Sol}}
\newcommand{\Hperv}{{}^\mathrm{p}\negmedspace H}

\newcommand{\Rb}{\mathbf{R}}

\newcommand{\et}{\text{\'et}}
\newcommand{\Fdepth}{F\negthickspace\op{--depth}}

\newcommand{\draft}[1]{}

\begin{document}

%
%

\title[]{Local cohomology multiplicities in terms of \'etale cohomology}
\author{Manuel Blickle}
\address{Universit\"at Essen, FB6 Mathematik, 45117 Essen,
Germany} \email{manuel.blickle@uni-essen.de} \urladdr{\url{www.mabli.org}}
\author{Raphael Bondu}
\address{Universit\"at Essen, FB6 Mathematik, 45117 Essen,
Germany} \email{raphael.bondu@uni-essen.de} 
\thanks{The authors were supported by the DFG Schwerpunkt \emph{Globale
Methoden in der komplexen Geometrie}.}
%
%

\maketitle

\begin{abstract}
In this paper we give an interpretation of the invariants $\lambda_{a,i}(A)$
introduced by Lyubeznik in \cite{Lyub.FinChar0} for a reasonably general class
of singularities. In positive characteristic it is the newly introduced class
of \emph{close to $F$--rational} varieties and the invariants are described in
terms of \'etale cohomology with $\ZZ/p\ZZ$--coefficients. This result presents
the first application of Emerton and Kisin's Riemann--Hilbert type
correspondence to local algebra. In fact our proof works in characteristic zero
as well so that we obtain generalizations of results on these invariants
which were previously obtained for isolated singularities by analytic
techniques.
\end{abstract}


%
%

\section{Introduction}
Let $(R,m)$ be a regular local ring of dimension $n$ and let $A=R/I$ be a
quotient of $R$. In \cite{Lyub.FinChar0} Lyubeznik introduces new
invariants $\lambda_{a,i}(A)$ (defined as the $a$th Bass number of
$H^{n-i}_I(R)$) and shows that if $A$ contains a field, these are independent of the presentation of $A$
as a quotient of a regular local ring. One can verify that
\[
    \lambda_{a,i}(A) = e(H^a_m(H^{n-i}_I(R)))
\]
where the multiplicity $e(\usc)$  can be described as follows: The main results of
\cite{Lyub.FinChar0,HuSha.LocCohom} state that the module $H^a_m(H^{n-i}_I(R))$
is injective. As it is supported at the maximal ideal it is isomorphic to a
finite direct sum of $e$ copies of the injective hull $E_{R/m} \cong H^n_m(R)$
of the residue field of $R$. This integer $e$ is the multiplicity.

Our main result is the following description of these invariants in the case that
$A$ has reasonable singularities.

\begin{theorem}\label{thm.main}
    Let $k$ be a field\footnote{As was pointed out by Brian Conrad there
is a problem (proof of Lemma \ref{lem.stupid}) with our argument (in characteristic
$p>0$) if $k$ is not separably closed. So ``$k$ separably closed'' has to be added
as an assumption in the theorem -- to be safe as an assumption throughout the paper.
\par
    However, in general one can reduce to this case. For $k'$ the separable
closure of $k$ one checks that $\lambda_{a,i}(A)=\lambda_{a,i}(A \tensor_k k')$ as
follows: Since tensoring with $k'$ over $k$ (which we shall denote by $(\usc)'$) is
faithfully flat, we have
\[
    H^{a}_{m'}(H^{n-i}_{I'}(R') \cong k' \tensor_k H^a_{m}H^{n-i}_I(R)
\]
and $E_{R/m} \cong E_{R'/m'}$ (using that $E_{R/m} \cong H^n_m(R)$) which implies
the claimed equality. Hence, in general, one is still able to reduce to the case
proven by our main theorem.} of characteristic $p>0$ and let $Y$
    be a $k$-variety of dimension $d$ which is close to $F$--rational outside
    the single (closed) point $x \in Y$. Let $A=\Oo_{Y,x}$. Then
\begin{enumerate}
    \item $\lambda_{0,i}(A) = \dim_{\ZZ/p\ZZ} H^i_{\{x\}}(Y_{\et},\ZZ/p\ZZ)$ for $ 1 \leq i \leq
    d-1$.
    \item $\lambda_{a,d}(A) = \dim_{\ZZ/p\ZZ} H^{d+1-a}_{\{x\}}(Y_{\et},\ZZ/p\ZZ)$ for $2 \leq a \leq
    d-1$ and $\lambda_{d,d}(A) = \dim_{\ZZ/p\ZZ} H^{1}_{\{x\}}(Y_{\et},\ZZ/p\ZZ)+1$.
    \item All other $\lambda_{a,i}(A)$ vanish.
\end{enumerate}
\end{theorem}
The notion \emph{close to $F$--rational} will be introduced and briefly
discussed in Section \ref{sec.Frat}. The name is chosen to indicate that
$F$--rational varieties are close to $F$--rational and thus so are smooth
varieties. In particular, the theorem applies in the case that $Y$ has an
isolated singularity at $x$. A key ingredient in our proof is that a close to
$F$--rational variety $Y \subseteq X$ ($X$ smooth) has the property that
$\Rb\Gamma_{[Y]}(\Oo_X)$ is isomorphic to $\Ll(Y,X)[d-n]$, the intersection
homology module (cf.\ Section \ref{sec.char0}).

The isolated singular case was motivated by the main result in
\cite{LopezSabbah} where Garc\'{i}a~López and Sabbah prove a topological
description of the invariants $\lambda_{a,i}(A)$ in the case that $A$ is the local ring of an isolated complex
singularity.\footnote{To be precise, they state part (2) in its Poincar\'e dual
form ($\lambda_{a,d}=\dim_\CC H^{d+a}(Y,\CC)$ for $2 \leq a \leq d$), see in
Remark \ref{rem.poinc} why this is the case and also why we prefer our
version.} Our methods lead to a generalization of their result,
replacing the assumption of isolated singularity with the significantly weaker requirement that $(Y-\{x\})$ is an \emph{intersection cohomology manifold}, \ie $\CC_{(Y-\{x\})}[d] \cong
IC^\bullet_{(Y-\{x\})}$.\footnote{This notion was recently introduced by Massey
\cite{Massey.IntCohom} and via the Riemann--Hilbert correspondence it is
clearly the condition corresponding to $\Ll(Y,X)\cong
\Rb\Gamma_{[Y]}(\Oo_X)[n-d]$}
\begin{theorem}\label{thm.char0}
    Let $Y$ be a complex variety of dimension $d$ and $x \in Y$ such that $(Y-\{x\})$
    is an intersection homology manifold. Let $A=\Oo_{Y,x}$. Then
    \begin{enumerate}
    \item $\lambda_{0,i}(A) = \dim_{\CC} H^i_{\{x\}}(Y,\CC)$ for $ 1 \leq i \leq
    d-1$.
    \item $\lambda_{a,d}(A) = \dim_{\CC} H^{d+1-a}_{\{x\}}(Y,\CC)$ for $2
    \leq a \leq d-1$ and $\lambda_{d,d}(A) = \dim_{\CC} H^{1}_{\{x\}}(Y,\CC)+1$.
    \item All other $\lambda_{a,i}(A)$ vanish.
\end{enumerate}
\end{theorem}
In the isolated singular case, statement (1) was already pointed out in
\cite{Lyub.FinChar0} to follow from a result of Ogus \cite[Theorem
2.3]{Ogus.LocCohm}. Observing the proof in \cite{LopezSabbah} we first note
that part (3) is independent of the characteristic whereas the other parts
distinctively use characteristic zero.

In order to obtain the full analog of Theorem \ref{thm.char0} in positive
characteristic we have to work somewhat harder. The proof given in
\cite{LopezSabbah} is our point of departure. They use the Riemann--Hilbert
correspondence and duality for holonomic $D$--modules. Our idea is to replace the
Riemann--Hilbert correspondence (\ie de Rham theory) with the correspondence
recently introduced by Emerton and Kisin \cite{EmKis.Fcrys} (\ie Artin--Schreier
theory). The main obstacle is that the categories involved in the Emerton--Kisin
correspondence do not have a duality, which was an essential part in the proof of
Garc\'{i}a~López and Sabbah. Thus our first task is to give a new proof of Theorem
\ref{thm.char0} which as its main feature avoids the use of duality. In this proof
we also show explicitly that part (1) and (2) are equivalent once part (3) is
established.

In Section \ref{sec.charp} we briefly recall the setup for the
Emerton--Kisin correspondence and show that this allows us to
carry our new characteristic zero proof over to positive
characteristic. Thus we obtain Theorem \ref{thm.main} postponing
the introduction and investigation of \emph{close to
$F$--rational} singularities to the last section.

\section{Duality free proof in characteristic
zero.}\label{sec.char0}\nopagebreak[4]
\subsection{A spectral sequence computation}\label{sec.spse}
We start with explaining that a vanishing
condition (slightly weaker than the one in part (3) of Theorem \ref{thm.char0})
for the $\lambda_{a,i}(A)$ implies part (3) and also that part (1) and (2) are
equivalent. This is done via a not so difficult spectral sequence argument. The
condition we impose is that
\begin{equation}\label{eq.vancond}
    \lambda_{a,i} = e(H^a_{[x]}H^{n-i}_{[Y]}(\Oo_X)) = 0
\end{equation}
for all pairs $(a,i)$ with $a>1$ and $i\neq d$. This is, for example, satisfied if $Y$
a complete intersection at $x$. We will then show that (\ref{eq.vancond})
implies that
\begin{gather*}
\lambda_{0,1}(A)+1 = \lambda_{d,d}(A) \quad \text{ and} \\
\lambda_{0,d-r+1}(A) = \lambda_{r,d}(A) \text{ for $2 \leq r \leq
d-1$}
\end{gather*}
and that all other $\lambda_{a,i}$ are zero. This clearly suffices to support
all our claims. Now consider the spectral sequence
\[
    E_2^{a,j} = H^a_{[x]}H^{j}_{[Y]}(\Oo_X) \Rightarrow
    H^{a+j}_{[x]}(\Oo_X).
\]
Since $\lambda_{a,i}=e(E_2^{a,n-i})$ the vanishing assumption
(\ref{eq.vancond}) yields that the only possibly nonzero entries
of the $E_2$ sheet of this spectral sequence are the ones
illustrated in the picture:
\[
\scriptsize{ \xymatrix@R0.1pc@C0.1pc{
    E_2^{0,n}\ar@{..}[dd] \\
    {} \\
    E_2^{0,n-d+2} \\
    E_2^{0,n-d+1}\ar[rrd] \\
    E_2^{0,n-d} & {E_2^{1,n-d}} & {E_2^{2,n-d}} & {E_2^{3,n-d}}\ar@{..}[rr] & &
    {E_2^{d-1,n-d}} & {E_2^{d,n-d}} \\
    E_2^{0,n-d-1}\ar@{..}[dd] \\
    {}\\
    E_2^{0,0}
}}
\]
Clearly, the only possibly nonzero arrow is the one indicated. We
now assume that $d \geq 2$ and leave the easy cases $d=1$ and
$d=0$ to the reader. Recall that by the Hartshorne--Lichtenbaum
vanishing theorem one has $H^{n}_{[Y]}(\Oo_X)=0$ and therefore
$E_2^{0,n}=H^0_{[x]}H^{n}_{[Y]}(\Oo_X)=0$ which just says that
$\lambda_{0,0}(A)=0$. Now we claim that for $r \geq 2$ the $E_r$
sheet of the spectral sequence has only the following (possibly)
nonzero entries
\[
\scriptsize{ \xymatrix@R0.1pc@C0.1pc{
    0 \\
    E_2^{0,n-1}\ar@{.}[dd] \\
    {}\\
    E_2^{0,n-d+r} \\
    E_2^{0,n-d+r-1} \ar[dddrrrr] \\
    {}\\
    {}\\
    0&&&& {E_2^{r,n-d}} & {E_2^{r+1,n-d}}\ar@{..}[rr] &&
    {E_2^{d-1,n-d}} & {E_2^{d,n-d}}
}}
\]
where the only nonzero arrow is the one indicated which yields an
isomorphism
\begin{equation}\label{eq.iso}
    H^0_{[x]}H^{n-d+r-1}_{[Y]}(\Oo_X) \cong H^r_{[x]}H^{n-d}_{[Y]}(\Oo_X)
\end{equation}
for $r<d$. In the case $r=d$ this only nonzero arrow fits in a
short exact sequence
\begin{equation}\label{eq.ses}
    0 \to H^0_{[x]}H^{n-1}_{[Y]}(\Oo_X) \to H^d_{[x]}H^{n-d}_{[Y]}(\Oo_X) \to
    H^n_{[x]}(\Oo_X) \to 0
\end{equation}
the right map being the edge map of the spectral sequence.\footnote{The
injectivity on the left is clear since $E^{0,n-1}_2$ must die in the limit. The
surjectivity on the right follows since $E_2^{0,0}=0$, thus the term in the
middle ($E^{d,d}_2$) is the only one that can contribute to the abutment term,
thus has to surject onto it.} All these claims simply follow from the
observation that in the limit of the spectral sequence the only surviving term
is $H^n_{[x]}(\Oo_X)$ and the fact that each sheet has only one nonzero arrow.
For $r<d$ the bottom left terms (the ones below the nonzero
arrow) must vanish since they do not contribute to the limit (the only term
that does is $E^{d,n-d}$) and since there are no nonzero arrows arriving at or
departing from any of them in the $r$th or any higher sheet. Similarly the only
nonzero arrow must be an isomorphism since otherwise there would be a surviving
term in the next sheet which is impossible as we just argued. Rephrasing these
observations in terms of the $\lambda_{a,i}(A)$ we obtain from
(\ref{eq.iso}) that
\[
    \lambda_{0,d-r+1}(A) = \lambda_{r,d}(A)
\]
for $2 \leq r \leq d-1$ and from (\ref{eq.ses}) that
\[
    \lambda_{0,1}(A)+1 = \lambda_{d,d}(A)
\]
where we used in the latter that $e(\usc)$ is additive on short exact sequences
and that $e(H^n_{[x]}(\Oo_X))=1$. That all other $\lambda_{a,i}$ vanish follows
already form the shape of the $E_2$--sheet.

\begin{remark}\label{rem.isoss}
The vanishing condition (\ref{eq.vancond}) is satisfied precisely when
$H^j_{[Y]}(\Oo_X)$ is supported at the point $x$ for $j \neq n-d$. This, in
turn, clearly holds whenever $Y$ has an isolated singularity at $x$ and is
smooth otherwise.
\end{remark}

\begin{remark}\label{rem.poinc}
In \cite{LopezSabbah} Garc\'{i}a~López and Sabbah prove the Poincar\'e
dual statement of part (2), namely that $\lambda_{a,d}=\dim_\CC H^{d+a}(Y,\CC)$
for $2 \leq a \leq d$.\footnote{In order to recover our part (2) of Theorem
\ref{thm.char0} one proceeds as in \cite[Remark 1]{LopezSabbah}) and uses
Poincar\'e duality for the link $L_{(Y,x)}$ of the singularity $(Y,x)$. The
link is a real orientable compact manifold of dimension $2d-1$. We have,
locally analytically around $x$ that
\[
\begin{split}
   H^{d+a-1}(Y-\{x\},\CC)  &\cong H^{d+a-1}(L_{(Y,x)},\CC) \\
                     &\cong H^{d-a}(L_{(Y,x)},\CC) \qquad \mbox{(Poincar\'e duality)} \\
                     &\cong H^{d-a}(Y-\{x\},\CC).
\end{split}
\]
At the same time $H^{d+a}_{\{x\}}(Y,\CC) \cong H^{d+a-1}(Y-\{x\},\CC)$ and
$H^{d-a+1}_{\{x\}}(Y,\CC) \cong H^{d-a}(Y-\{x\},\CC)$ for $a \neq d$. When
$a=d$ one has $\dim_\CC H^{1}_{\{x\}}(Y,\CC) = \dim_{\CC} H^{0}(Y-\{x\},\CC)+1$
and the claim follows.} The reason for this lies in their computation of
$\lambda_{a,d}(A)$ which uses duality for holonomic $\Dd$--modules which under
Riemann--Hilbert corresponds to Poincar\'e duality, in that special case. Thus
they obtain the equivalence of part (1) and part (2) as a consequence of
Poincar\'e duality. Our observation though shows that this equivalence follows
from the structure of the invariants $\lambda_{a,i}$ and the use of Poincar\'e duality
can be avoided.
\end{remark}

\subsection{Preparatory lemmata and proof of Theorem \ref{thm.char0}}
We start with some (probably well known) facts which will naturally lead to the
proof of Theorem \ref{thm.char0}.
\begin{lemma}\label{lem.SolH}
    Let $X$ be a smooth $\CC$--variety of dimension $n$ and let $k:x \into X$ be
    the inclusion of a point. Let $\Mm$ be a holonomic $\Dd_X$--module, then
    \[
        \Sol(H^a_{[x]}(\Mm)) \cong k_{!}k^{-1}H^{-a}(\Sol \Mm).
    \]
\end{lemma}
\begin{proof} By definition of the symbols involved ($\Hperv$ denotes perverse
cohomology, $\Sol(\usc) \defeq \Rb \Hom_{\Dd_X}(\usc,\Oo_X)[n])$
we have
\[
    \Sol(H^a_{[x]}(\Mm)) \cong \Sol H^a(\Rb\Gamma_{[x]}\Mm) \cong
    \Hperv^{-a}(\Sol(\Rb\Gamma_{[x]}\Mm)).
\]
One has that $\Sol \circ \Rb\Gamma_{[x]} \cong k_!k^{-1} \circ \Sol$.\footnote{This
follows from the fact that for any closed embedding $k: Z \to X$ and complex of
$\Dd_X$--modules $\Mm^\bullet$ one has the triangle
\[
    \Rb\Gamma_{[Z]}\Mm^\bullet \to \Mm^\bullet \to \Rb j_* j^*
    \Mm^\bullet \to[+1]
\]
where $j$ denotes the open inclusion $X-Z \subseteq X$. Let us denote
$\Sol(\Mm^\bullet)$ by $\Ll^\bullet$ and apply $\Sol$ to this triangle. Using
the compatibility of $\Sol$ with the six operations, in particular $\Sol \circ
j_* \cong j_! \circ \Sol$ and $\Sol \circ j^{-1} \cong j^! \circ \Sol$, we obtain the
following triangle.
\[
    \Rb j_!j^!\Ll^\bullet \to \Ll^\bullet \to \Sol(\Rb\Gamma_{[Z]}\Mm^\bullet) \to[+1]
\]
Comparing with the standard triangle  $\Rb
j_!j^!\Ll^\bullet \to \Ll^\bullet \to \Rb k_!k^{-1} \Ll^\bullet
\to[+1]$ \cite[Triangle
2.6.33]{KashiwaraSchapira.SheavesManifolds} one obtains the claim. } Thus $\Hperv^{-a}(\Sol(\Rb\Gamma_{[x]}\Mm))$ is
equal to
\[
    \Hperv^{-a}(k_!k^{-1}(\Sol(\Mm))) \cong k_! \Hperv^{-a}
    (k^{-1}(\Sol(\Mm)))
\]
where the isomorphism holds since $k_!$ is $t$--exact by \cite[Lemma
III.4.1]{KielWeiss.WeilConj}. After pullback along $k$ we are on the point $x$
on which perverse cohomology is the same as ordinary cohomology. Thus we may
replace $\Hperv^{-a}$ by $H^{-a}$. Using that $k^{-1}$ is exact we get
\[
    k_! \Hperv^{-a} (k^{-1}(\Sol(\Mm))) \cong k_! H^{-a}(k^{-1}(\Sol(\Mm))) \cong
    k_!k^{-1} H^{-a}(\Sol(\Mm))
\]
as claimed.
\end{proof}
One of our tools is the intersection homology $\Dd_X$--module
$\Ll(Y,X)$ of Brylinski and Kashiwara \cite{Bry.Kash}. It is the
middle extension
\[
    \Ll(Y,X) \cong \tilde{j}_{!*}H^{n-d}_{[Y-\op{Sing} Y]}(\Oo_{(X-\op{Sing} Y)})
\]
where $\tilde{j}$ denotes the open inclusion $(X-\op{Sing} Y) \subseteq X$. Its
characterizing property is that it is the smallest $\Dd_X$--submodule of
$H^{n-d}_{[Y]}(\Oo_X)$ which agrees with $H^{n-d}_{[Y]}(\Oo_X)$ away from the
singular locus of $Y$. Thus in particular if $Y$ is smooth then
$\Ll(Y,X) = H^{n-d}_{[Y]}(\Oo_X)$.
\begin{lemma}\label{lem.intHom}
    Let $X$ be a smooth $k$--variety of dimension $n$, let $i:Y \into X$ be
    a closed subvariety of dimension $d$ and assume that for $x
    \in Y$ one has
    $\Ll(Y,X)|_{(X-\{x\})}\cong H^{n-d}_{[Y]}(\Oo_{X})|_{(X-\{x\})}$. Then
    \[
        H^a_{[x]}(H^{n-d}_{[Y]}(\Oo_X)) \cong H^a_{[x]}(\Ll(Y,X))
    \]
    for $a \geq 2$.
\end{lemma}
\begin{proof}
By assumption one has the short exact sequence
\[
    0 \to \Ll(Y,X) \to H^{n-d}_{[Y]}(\Oo_X) \to \Cc \to 0
\]
whose cokernel $\Cc$ has support in the point $x$. Thus
$H^a_{[x]}(\Cc)=0$ for $a\geq 1$. By the long exact sequence of
$H^\bullet_{[x]}(\usc)$ applied to this short exact sequence the
claim of the lemma follows.
\end{proof}

\begin{lemma}\label{lem.Soli!}
    Let $X$ be a smooth $k$--variety of dimension $n$ and let $Y \subseteq X$ be
    a closed subvariety of dimension $d$. Assume that for $x \in
    Y$ one has
    $\Ll(Y,X)|_{(X-\{x\})}\cong\Rb\Gamma_{[Y]}^{n-d}(\Oo_X)|_{(X-\{x\})}[n-d]$.
    Then
    \[
        \Sol(\Ll(Y,X)) \cong i_!j_{!*}\CC_{(Y-\{x\})}[d]
    \]
    where $j$ is the inclusion of $(Y-\{x\})\into Y$.
\end{lemma}
\begin{proof}
Let us fix the notation $(X-\op{Sing} Y) \to[j''] (X-\{x\})\to[j'] X$.
Then by definition of middle extension and our assumption one has
\[
\begin{split}
    \Ll(Y,X) &\cong j'_{!*}j''_{!*}H^{n-d}_{[Y-\op{Sing} Y]}(\Oo_{(X-\op{Sing}
    Y)})\\
             &\cong j'_{!*}\Ll(Y-\{x\},X-\{x\}) \\
             &\cong j'_{!*}\Rb\Gamma_{[Y-\{x\}]}(\Oo_{(X-\{x\})})[n-d].
\end{split}
\]
Denoting the inclusion $(Y-\{x\}) \into (X-\{x\})$ by $i'$ we have (see footnote~5)
\[
\begin{split}
    \Sol(\Rb\Gamma_{[Y-\{x\}]}(\Oo_{(X-\{x\})})[n-d])
        &\cong i'_!{i'}^{-1}\Sol(\Oo_{(X-\{x\})})[d-n]\\
        &\cong i'_!{i'}^{-1}\CC_{(X-\{x\})}[n][d-n] = i'_! \CC_{(Y-\{x\})}[d]
\end{split}
\]
where we also used that $\Sol(\Oo_{(X-\{x\})})\cong \CC_{(X-\{x\})}[n]$. Now
finish the proof with the following chain of equalities
\[
\begin{split}
    \Sol(\Ll(Y,X)) &\cong j'_{!*}\Sol(\Rb\Gamma_{[Y-\{x\}]}(\Oo_{(X-\{x\})})[n-d])\\
    &\cong j'_{!*}i'_{!}\CC_{(Y-\{x\})}[d] \\
    &\cong i_!j_{!*}\CC_{(Y-\{x\})}[d]
\end{split}
\]
the last of which follows from the fact that for a closed
immersion $i_!\cong i_{!*}$ and thus the $j$ and $i$ can be exchanged
as we have done.
\end{proof}
\begin{remark}
Granted, the assumption on $\Ll(Y,X)$ of the preceding two lemmata seems
somewhat random. In characteristic zero (say over $\CC$) they are equivalent
via the Riemann--Hilbert correspondence to $(Y-\{x\})$ being an intersection
cohomology manifold, see \cite{Massey.IntCohom}. In positive characteristic our
notion of close to $F$--rational of the final section relates it to previous
work on singularities, such as tight closure theory and the notion of
$F$--depth as in \cite{HaSp}.

Also note that if $Y$ has an isolated singularity at $x$ then the assumptions
are (trivially) satisfied since in this case one has
$\Ll(Y,X)|_{(Y-\{x\})}\cong H^{n-d}_{(Y-\{x\})}(\Oo_{(X-\{x\})}) \cong \Rb\Gamma_{[Y-\{x\}]}(\Oo_{(X-\{x\})})[n-d]$.
\end{remark}

\begin{proof}[Proof of Theorem \ref{thm.char0}]
By assumption $(Y-\{x\})$ is an intersection homology manifold which in particular
implies by Remark \ref{rem.isoss} that part (3) holds and part (1) and (2) are
equivalent. Thus it is enough to show, say, part (2):

As we pointed out in the introduction the $\Dd_X$--module
$H^a_{[x]}H^i_{[Y]}(\Oo_X)$ is isomorphic to a finite direct sum of
$\lambda_{a,i}(A)$ many copies of $H^n_{[x]}(\Oo_X)$, the injective hull of the
residue field at $x$. By Lemma \ref{lem.SolH} together with $\Sol(\Oo_X) =
\CC_X[n]$ one has
\[
    \Sol(H^n_{[x]}(\Oo_X)) \cong k_!k^{-1}H^{-n}(\CC_X[n]) \cong k_! \CC_x
\]
where we recall that $k$ was just the inclusion of $x \into X$.
Therefore $\lambda_{a,i}(A)$ is just the dimension of the fiber at
$x$ of $\Sol(H^a_{[x]}H^i_{[Y]}(\Oo_X))$. Thus, for $a \geq 2$ we
can compute
\[
\begin{split}
    \lambda_{a,d}(A) &= e(H^a_{[x]}(H^{n-d}_{[Y]}(\Oo_X))) \\
                     &= \dim_{\CC} \left(\Sol(H^a_{[x]}(H^{n-d}_{[Y]}(\Oo_X)))\right)_x \\
                     &= \dim_{\CC} \left(\Sol(H^a_{[x]}(\Ll(Y,X)))\right)_x  \qquad\text{ (by Lemma \ref{lem.intHom})}\\
                     &= \dim_{\CC} \left(k_!k^{-1} H^{-a}(\Sol
                     \Ll(Y,X))\right)_x \qquad \text{ (by Lemma
                     \ref{lem.SolH})} \\
                     &= \dim_{\CC} \left(H^{-a}(i_!j_{!*}\CC_{(Y-\{x\})}[d])
                     \right)_x \qquad \text{ (by Lemma \ref{lem.Soli!})} \\
                     &= \dim_{\CC} \left(H^{-a}(j_{!*}\CC_{(Y-\{x\})}[d])
                     \right)_x
\end{split}
\]
where $i$ is the inclusion $Y \into X$ and $j$ denotes the inclusion $(Y-\{x\})
\subseteq Y$. Since $j$ is just the inclusion of the complement of a single
point it follows that
\[
    j_{!*}\CC_{(Y-\{x\})}[d] \cong \tau_{\leq {d-1}} \Rb j_*\CC_{(Y-\{x\})}[d]
\]
by \cite[V.2.2 (2)]{Borel.Icohom}. By definition of Deligne's truncation
$\tau_{\leq {d-1}}$ one has for $a \geq 1$
\[
    \left(H^{d-a}(j_{!*}\CC_{(Y-\{x\})})\right)_x \cong \left(\Rb^{d-a}
    j_*\CC_{(Y-\{x\})}\right)_x.
\]
Applying the following Lemma \ref{lem.stupid} we get for $2 \leq a \leq d-1$
that
\[
    \lambda_{a,d}(A)= \dim_{\CC} \left(\Rb^{d-a}j_*\CC_{(Y-\{x\})}\right)_x =
    \dim_{\CC} H^{d-a+1}_{\{x\}}(Y,\CC)
\]
and (for $a=d$) that $\lambda_{d,d}(A)= H^{1}_{\{x\}}(Y,\CC)+1$.
\end{proof}

\begin{lemma}\label{lem.stupid}
    Let $Y$ be a variety and let $x \in Y$ be a closed point and $C$ be a constant
    sheaf on $Y$. Then
    \[
        (\Rb^ij_*j^{-1}C)_x \cong H^{i+1}_{\{x\}}(Y,C)
    \]
    for $i \geq 1$ and (for $i=0$) one has the short exact sequence
    \[
        0 \to C_x \to (\Rb^0j_*j^{-1}C)_x \to H^1_{\{x\}}(Y,C) \to 0.
    \]
\end{lemma}
\begin{proof}
    For the open inclusion $j:(Y-\{x\}) \into Y$ consider the triangle
    \[
        \Rb\Gamma_{\{x\}} C \to C \to \Rb j_*j^{-1} C \to[+1]
    \]
    and take its fiber at the point $x$ to obtain the following triangle:
    \[
        (\Rb\Gamma_{\{x\}} C)_x \to C_x \to (\Rb j_*j^{-1} C)_x \to[+1]
    \]
    Since $H^i(C_x)=0$ for $i>0$ and $H^{0}_{\{x\}}(C)=0$ (since $C$ is a
    constant sheaf) the long exact sequence  of cohomology for this triangle yields
    \[
        0 \to C_x \to (\Rb^0j_*j^{-1}C)_x \to (H^1_{\{x\}}(Y,C))_x \to 0
    \]
    and for $i \geq 1$
    \[
        (R^i j_*j^{-1}C)_x\cong (H^{i+1}_{\{x\}}(C))_x .
    \]
    But clearly since $H^{i+1}_{\{x\}}(C)$ is supported on $x$ we have
    $(H^{i+1}_{\{x\}}(C))_x  \cong H^{i+1}_{\{x\}}(Y,C)$ which finishes the proof.
\end{proof}

\section{The case of positive characteristic.}\label{sec.charp}
We very briefly recall the setup of the correspondence of Emerton and Kisin and
point out the relevant facts which will make clear that the  in
characteristic zero also works in positive characteristic.

\subsection{Emerton--Kisin correspondence}
Let $k$ be a field of positive characteristic $p$ and let $X$ be a smooth
$k$--variety. In \cite{EmKis.Fcrys} Emerton and Kisin establish an
anti--equivalence (on the level of derived categories) between constructible
${\ZZ/p\ZZ}$--sheaves on $X_\et$ on one hand and locally finitely generated
unit $\Oo_{F,X}$--modules on the other. Their construction closely models the
Riemann--Hilbert correspondence and underlies the same formalism -- except that
there is no duality available on either side of the correspondence. This leads
to the defect that their anti--equivalence is compatible with only three of
Grothendieck's \emph{six operations}, namely with analogs of $f^!$, $f_*$
(denoted $f_+$ in \cite{EmKis.Fcrys}) and $\overset{\LL}{\tensor}_{\Oo_{F,X}}$
on the $\Oo_{F,X}$--module side, which correspond to $f^*$, $f_!$ and
$\overset{\LL}{\tensor}_{{\ZZ/p\ZZ}}$ on the constructible \'etale ${\ZZ/p\ZZ}$
side.

We will recall the definition of $\Oo_{F,X}$--module and point out
that the local cohomology modules are locally finitely generated
as such, so that the formalism of Emerton--Kisin can be applied to
our study of the numbers $\lambda_{a,i}$. For a nice introduction
see \cite{EmKis.FcrysIntro}; or \cite{EmKis.Fcrys} for the most
general theory.
\begin{definition}
A quasi--coherent $\Oo_X$--module $\Mm$ together with an $\Oo_X$--linear map
\[
    \theta: F^*\Mm \to \Mm
\]
is called an \emph{$\Oo_{F,X}$--module}. Here $F$ denotes the Frobenius
morphism on $X$. If $\theta$ is an isomorphism, then $(\Mm,\theta)$ is called
\emph{unit}.
\end{definition}
Locally, if $X=\Spec R$, an $\Oo_{X,F}$--module is nothing but a module $M$
over the non-commutative ring
\[
    R[F] \defeq \frac{R \langle F \rangle}{(r^pF-Fr\;|\; r \in R)}.
\]
Such $R[F]$--module is called \emph{finitely generated} if it is just that:
finitely generated as an $R[F]$--module. Thus we have the notion of locally
finitely generated for $\Oo_{F,X}$--modules.

One of the key results of the theory (which was essentially proved by Lyubeznik
in \cite{Lyub}) is that the category of locally finitely generated unit
$R[F]$--modules is abelian, and that every such $\Mm$ has finite length in that
category \cite[Theorem 3.2]{Lyub}.

\begin{example}
    Let $X=\Spec R$ be affine. Then, abusing the identification of
    $\Oo_X$--modules and $R$--modules one can write $F^*M =
    R^{(1)} \tensor_R M$ where $R^{(1)}$ is the $R$--$R$--bimodule
    with the usual left structure and the right structure via the
    Frobenius map. Thus one sees immediately that the natural map
    \[
        F^*R = R^{(1)} \tensor_R R \to R
    \]
    sending $a \tensor r$ to $ar^p$ is an isomorphism, showing
    that $R$ is a fg (finitely generated) unit $R[F]$--module.

    Let $g \in R$ be an element and consider the localization
    $R_g$. The natural map
    \[
        F^*R_g =R^{(1)} \tensor_R R_g \to R_g
    \]
    has an inverse given by sending $r/t$ to $rt^{p-1} \tensor
    1/t$.  $R_g$ is generated as an $R[F]$--module by
    $1/g$. Thus again $R_g$ is a fg unit $R[F]$--module.

    Since local cohomology modules $H^i_I(R)$ for $I$ an ideal of $R$
    can be computed via \v{C}ech resolutions,
    whose entries are localizations of the type $R_g$, the aforementioned
    result that the category of fg unit $R[F]$--modules is abelian implies
    that local cohomology modules are fg unit.
\end{example}

These examples are a special instance of more general results showing that the
\emph{cohomology with supports} functors are defined in the category of locally
finitely generated unit $\Oo_{F,X}$--modules \cite[Proposition
5.11.5]{EmKis.Fcrys}. If $\Mm^\bullet$ is a bounded complex of such modules
then so is $\Rb\Gamma_{[Y]}\Mm^\bullet$ for $Y$ a closed subvariety of $X$ and
one has the usual triangle
\[
    \Rb\Gamma_{[Y]}\Mm \to \Mm \to j_+j^!\Mm \to[+1]
\]
where $j: X-Y \to X$ denotes the open inclusion.

The correspondence of Emerton and Kisin is between the derived category of
bounded complexes of $\Oo_{F,X}$ modules whose cohomology is locally finitely
generated unit,
\[
\xymatrix{
    {D^b_{lfgu}(\Oo_{F,X})} \ar@/^/^{\Sol}_{\sim}[rr] &&\ar@/^/[ll] {D^b_c(X_{\et},{\ZZ/p\ZZ})} }
\]
and the derived category of bounded complexes of ${\ZZ/p\ZZ}$ sheaves with
constructible cohomology on $X_{\et}$. Furthermore, they define functors $f^!$,
$f_+$ and $\usc\overset{\LL}{\tensor}_{\Oo_{F,X}}\usc$. They are not the same
as (though closely related to) the functors of Grothendieck-Serre duality.

The canonical $t$--structure on the left induces via the anti--equivalence an
exotic $t$--structure on $D^b_c(X_{\et},{\ZZ/p\ZZ})$, which in turn is just the
$t$--structure for the middle perversity as described by Gabber
\cite{Gabb.tStruc}. Thus one obtains a notion of \emph{perverse sheaves} and
thus of perverse cohomology.

\subsection{Intermediate extensions}
There is a theory of intermediate extensions. If $j: U \to X$ is a locally
closed immersion of smooth $k$--schemes and $\Mm$ a lfgu $\Oo_{F,U}$--module,
then its intermediate extension $j_{!+}\Mm$ is defined as the smallest
submodule $\Mm' \subseteq H^0(j_+\Mm)$ such that $j^!\Mm'=\Mm$. Furthermore,
        \[
            \Sol(j_{!+}(\Mm)) \cong j_{!*}\Sol(\Mm) \defeq
            \Image(\Hperv^0(j_!\Sol(\Mm)) \to \Hperv^0(j_*\Sol(\Mm)))
        \]
so that the intermediate extension is compatible with the correspondence
\cite[Section 4.3]{EmKis.FcrysIntro}. We will only apply this to
$\Mm=H^{n-d}_{[Y-\op{Sing} Y]}(\Oo_{(X-\op{Sing} Y)})$ for $Y$ a closed subset
of $X$ and $j$ the open inclusion $(X-\op{Sing} Y) \subseteq X$. In this case
we get
\[
    \Ll(Y,X) \defeq j_{!+}H^{n-d}_{[Y-\op{Sing} Y]}(\Oo_{(X-\op{Sing} Y)}) \subseteq H^{n-d}_{[Y]}(\Oo_X)
\]
as its unique simple submodule. This important special case was already
obtained in \cite{Manuel.int}. The key point in obtaining these results is the
aforementioned fact that lfgu $\Oo_{F,X}$--modules have finite length.

The following proposition lists the properties of the theory which are needed
to be able to transfer the proof of Theorem \ref{thm.char0} to positive characteristic.
\begin{proposition}\label{prop.properties}
\begin{enumerate}
   \item $\Sol(\Oo_X) \cong {\ZZ/p\ZZ}[n]$ where $n$ is the dimension of $X$.
    \item For a closed immersion of smooth $k$--schemes $k: Y \to X$ one has $\Sol \circ
\Rb\Gamma_{[Y]} \cong k_!k^{-1} \circ \Sol$.
    \item Let $k: Y \to X$ be a closed immersion of smooth schemes. Then $k_!$ is $t$--exact.

\end{enumerate}
\end{proposition}
\begin{proof}
Part (1) is just Example 9.3.1 in \cite{EmKis.Fcrys}.

For part (2) note that $\Rb\Gamma_{[Y]}\Mm$ is defined via the triangle
\[
    \Rb\Gamma_{[Y]}\Mm \to \Mm \to j_+j^!\Mm \to[+1]
\]
with $j: X-Y \to X$ denoting the open inclusion. Applying $\Sol$ and using the
fact that $\Sol$ interchanges $j_+$ with $j_!$ and $j^!$ with $j^*$ by
\cite[Proposition 9.3, Proposition 9.5]{EmKis.Fcrys} we compare with the
triangle
\[
    j_!j^*\Sol(\Mm) \to \Sol(\Mm) \to k_!k^{-1}\Sol(\Mm) \to[+1]
\]
in $D^b_c(X_\et, {\ZZ/p\ZZ})$ to obtain the result.

Part (3) can be checked by hand (using Gabbers definition of the $t$--structure
in \cite{Gabb.tStruc}), but also follows via the correspondence from the fact
that $k_+$ is exact by \cite[Remark 3.4.1]{EmKis.Fcrys}.
\end{proof}

\begin{proof}[Proof of Theorem \ref{thm.main}]
The assumption of close to $F$--rational implies by Proposition
\ref{prop.clFrat} that $\Ll(Y,X)|_{(X-\{x\})} \cong
\Rb\Gamma_{Y}(\Oo_{X})|_{(X-\{x\})}[n-d]$. This means in particular that
$\Ll(Y,X)|_{(X-\{x\})} \cong H^{n-d}_{Y}(\Oo_{X})|_{(X-\{x\})}$ and $H^{i}_{Y}(\Oo_{X})$ is
supported in $x$ for $i \neq n-d$. Thus the vanishing condition
(\ref{eq.vancond}) is satisfied and therefore (by Section \ref{sec.spse}) part
(3) holds and part (1) and (2) are equivalent. Again we prove part (2) to
finish the argument.

This is done by following the arguments in the preceding section step by step,
working on the \'etale site and replacing $\CC$ by $\ZZ/p\ZZ$ whenever
appropriate. Here are some remarks on this task which finishes the proof.
\begin{enumerate}
\item For Lemma \ref{lem.SolH} one uses that $\Rb\Gamma_{[x]}$ commutes with
    $\Sol$ in the way claimed. Furthermore we use that $k_!$ is $t$--exact. This is
    Proposition \ref{prop.properties} part (2) and (3).
\item As pointed out at the beginning of the proof the assumptions of Lemma
    \ref{lem.intHom} and Lemma \ref{lem.SolH} are satisfied. For Lemma \ref{lem.intHom}
    literally the same argument holds after the existence of the middle extension $\Ll(Y,X)$ is
    established as discussed above. The same remark applies to Lemma \ref{lem.Soli!}.
\item In the actual proof one should use that $\Sol(\Oo_X) \cong \ZZ/p\ZZ[n]$ and the
    discussed properties of middle extension, in particular its compatibility with
    $\Sol$.
\item Lemma \ref{lem.stupid} is even stated for general coefficients and the
    argument is valid for any $k$--variety with the \'etale topology. The only caveat is that we
    implicitly used excision in the last part; an \'etale version of which
    can be found in \cite[Chapter 3, Proposition 1.27]{Milne}, for
    example.\footnote{As Brian Conrad pointed out to me, this use of excision is more
    subtle than we indicate here if $k$ is not separably closed. Thus to be correct we
    assume $k$ separably closed.}
\end{enumerate}
\end{proof}

\begin{remark}
In positive characteristic a more direct proof of our main result is possible.
One observes that the numbers $\lambda_{0,i}(A)$ for $i=1 \ldots d-1$ can be
interpreted in terms of the action of the Frobenius on $H^i_{m}(A)$. Namely if
$A=R/I$ the local cohomology module $H^{n-i}_I(R)$ is obtained from the local
cohomology module $H^i_m(A)$ via a certain functor $\Hh_{R,A}$ (introduced and
studied in \cite[Section 4, Example 4.8]{Lyub}). This functorial relationship
        \[
            H^{n-i}_I(R) \cong \Hh_{R,A}(H^i_m(A))
        \]
implies that $\lambda_{0,i}(A)=e(H^{n-i}_I(R))$ (which is called the
\emph{corank} of $H^{n-i}_I(R)$ in \cite{Lyub}) is equal to the dimension of
the Frobenius stable part of $H^i_m(A)$, by \cite[Proposition 4.10]{Lyub}

    Finally, in \cite[Theorem 2.5]{HaSp} the latter is determined to be
    equal to $\dim_k H^i_{\{x\}}(Y_\et,\ZZ/p\ZZ)$ as required.
\end{remark}

\section{Close to $F$--rational singularities.}\label{sec.Frat}
We finish this note with a brief discussion of a new class of singularities,
called close to $F$--rational.
\begin{definition}
Let $(A,m)$ be a local noetherian ring of dimension $d$. Let
$H^\bullet_m(A)=\dirsum H^i_m(A)$ be the local cohomolgy with support in $m$.
Then $A$ is called \emph{close to $F$--rational} if and only if
\[
    H^*_m(A)/0^F
\]
is simple as an $A[F]$--module, where $0^F$ denotes all the elements of
$H^*_m(A)$ which are annihilated by a power of the Frobenius.

If $Y$ is a noetherian scheme then $Y$ is called \emph{close to $F$--rational}
if for all closed points $y \in Y$ the local ring $\Oo_{Y,y}$ is close to
$F$--rational.
\end{definition}
Recall that $F$--rationality of $A$ is equivalent to $H^*_m(A)$ being simple as
an $A[F]$--module (at least if $A$ is excellent). This implies that an
$F$--rational ring is close to $F$--rational. The obstruction to
$F$--rationality is the tight closure of zero $0^*$ in $H^*_m(A)$ (see
\cite{Hune.tight} for relevant notions from the theory of tight closure). Close
to $F$-rational just means that this obstruction is, if not zero
($F$--rational) it is at least $F$--nilpotent. Since one always has that
$H^d_m(A) \neq 0^*$ \footnote{For the versed in tight closure theory: existence
of test elements is responsible for this, see for example \cite{Hune.tight}.}
it follows that $A$ is almost $F$--rational if and only if
$H^d_m(A)/0_{H^d_m(A)}^F$ is $A[F]$--simple and $H^i_m(A)$ is $F$--nilpotent for
$i \neq d$.

The following characterization of close to $F$--rational singularities is the
key point of our investigation.
\begin{proposition}\label{prop.clFrat}
    Let $Y$ be a subvariety of dimension $d$ embedded in $X$ which is a smooth
    $k$--variety of dimension $n$ ($\op{char} k = p >0$). Then $Y$ is close to
    $F$--rational if and only if
    \[
        \Ll(Y,X) \cong \Rb\Gamma_{[Y]}(\Oo_X)[n-d]
    \]
    where $\Ll(Y,X)$ denotes the unique simple unit $R[F]$ submodule of
    $H^c_{[Y]}(\Oo_X)$.
\end{proposition}
\begin{proof}
    This is a slight extension of the main result in \cite{Manuel.int}. Since
    by definition, close to $F$--rational is checked locally, we have to verify
    that for every point $y \in Y$ the local ring $(A,m) = \Oo_{Y,y}$ is
    close to $F$--rational if and only if $\Ll(A,R) = \Rb\Gamma_{I}(R)[n-d]$ where
    $R=\Oo_{X,x}$ such that $A=R/I$. For this reduction we used that $\Ll(Y,X)$
    and $\Rb\Gamma_{[Y]}\Oo_X$ localize.

    In this situation \cite[Theorem 4.9]{Manuel.int} states that
    $\Ll(A,R)=H^{n-d}_I(R)$ if and only if $0^*=0^F$
    holds in $H^d_m(A)$. This latter condition is equivalent to $H^d_m(A)/0^F$
    being $A[F]$--simple since $0^*$ is the maximal proper $R[F]$--submodule of
    $H^d_m(A)$. It remains to point out that $H^{n-i}_I(R)$ is zero
    if and only if $H^{i}_m(A)$ is $F$--nilpotent. This is because, in the
    notation of \cite[Example 4.8]{Lyub} we have
    \[
        H^{n-i}_I(R) \cong \Hh_{R,A}(H^i_m(A)).
    \]
    By \cite[Section 4]{Lyub} one has $\Hh_{R,A}(\Mm)=0$ if and only if the
    $A[F]$--module $\Mm$ is $F$--nilpotent. It follows that $H^{n-i}_I(R)=0$ if
    and only $H^i_m(A)$ is $F$--nilpotent.
\end{proof}

\begin{remark}
Close to $F$--rational singularities are related to the notion of
$\Fdepth$ as introduced by Hartshorne and Speiser in \cite[page
60]{HaSp}. One can verify that if $A$ is close to $F$--rational
then $\Fdepth A = \dim A$. This notion of $\Fdepth$ is shown to be
equal to the \'etale $\ZZ/p\ZZ$--depth.

Thus the notion of (close to) $F$--rational singularities yields a reasonable
description of the class of varieties $Y \subseteq X$ for which $\Ll(Y,X) \cong
\Rb\Gamma_{[Y]}(\Oo_X)[n-d]$ and consequently Theorem \ref{thm.main} holds.
\end{remark}

\begin{remark}
    To see that close to $F$--rational does not imply $F$--rational one can consider the example of $A=k[x,y,z]/(x^4+y^4+z^4)$. This is not $F$--rational but in \cite[Example 5.28]{Manuel.PhD} the first author shows that it is close to $F$--rational precisely if the characteristic $p$ of $k$ is congruent to $3$ mod $4$.
\end{remark}

\begin{remark}
In a recent preprint \cite{Massey.IntCohom} Massey introduces and studies the
notion of \emph{intersection homomolgy manifold} in characteristic zero. As we
pointed out before, this notion means precisely that
\[
    \Ll(Y,X) \cong \Rb\Gamma_{[Y]}(\Oo_X)[n-d]
\]
whenever $Y \subseteq X$ is embedded into a smooth variety $X$. He gives
several alternative characterisations, particularly for $Y$ a complete
intersection.

Also in the complete intersection case Torrelli gives in
\cite{TorrelliT.IntHomBernstein} a characterization of $\Ll(Y,X) \cong
\Rb\Gamma_{[Y]}(\Oo_X)[n-d]$ in terms of the Bernstein polynomials. In the case
that $Y=(f=0)$ is a hypersurface his condition is easily phrased: The reduced
Bernstein polynomial (that is divide the usual Bernstein polynomial by $(x+1)$)
of $f$  has no integral root $\leq -1$.
\end{remark}

\providecommand{\bysame}{\leavevmode\hbox to3em{\hrulefill}\thinspace}
\providecommand{\MR}{\relax\ifhmode\unskip\space\fi MR }
\providecommand{\MRhref}[2]{%
  \href{http://www.ams.org/mathscinet-getitem?mr=#1}{#2}
} \providecommand{\href}[2]{#2}

%
%

\end{document}